# A Graceful Algebraic Function Labelling of Rooted Symmetric Trees

RAFAEL I. ROFA

**ABSTRACT.** Let $T = (V,E)$ be a tree with vertex set V and edge set E. A graceful labelling f of T is an injective function $f : V \to \{0, 1, \ldots, |E|\}$ such that if edge uv is assigned the label $g(uv)=|f(u)-f(v)|$ then the function $g : E \to \{1, \ldots, |E|\}$ is also injective (that is all edge labels are distinct). A rooted symmetric tree is a tree in which all vertices at the same level from root vertex have the same degree. It has been known since 1979 that rooted symmetric trees are graceful. However, the proofs that have been presented for this fact are either indirect inductive proofs or algorithmic descriptive proofs showing the many separate steps involved in labelling the vertices. Given a rooted symmetric tree, we find a graceful labelling f of T in the form of a direct algebraic function that algebraically *maps* each vertex of T to a unique label. Interestingly, f turns out to be a generalisation of the way a path is 'canonically' gracefully labelled as *algorithmically* described by A. Rosa who is a pioneer in this field of research. The algebraic function that defines the graceful labelling is used to show that a class of rooted symmetric trees that contains the class of binomial trees has weakly α-labelling and it can provide a concise practical, computational way of producing graceful labellings of 'large' rooted symmetric trees in, relatively, minimal time.

## 1. Introduction

A graceful labelling of a tree $T = (V, E)$ with vertex set V and edge set E, is an injective function $f : V \to \{0, 1, \ldots, |E|\}$ such that if edge uv is assigned the label $g(uv)=|f(u)-f(v)|$ then the function $g : E \to \{1, \ldots, |E|\}$ is also injective (that is all edge labels are distinct) . Several graph labelling schemes were first introduced by A. Rosa (Rosa, 1967) as a means of attacking the Ringel Conjecture that the complete graph $K_{2n+1}$ can be decomposed into 2n+1 isomorphic trees each on n edges (Ringel, 1963). Although A. Rosa introduced four different labelling schemes of graphs that can lead to the proof of the Ringel conjecture, nevertheless β-valuations, or graceful labelings, as were subsequently called by Golomb (1972), caught the attention and efforts of the largest part of researchers resulting in a vast reservoir of papers about the subject. For the purpose of this paper we also define an α-labelling of a tree T to be a graceful labelling with the additional property that there is a positive integer k such that for any edge xy either $f(x) \leq k < f(y)$ or $f(y) \leq k < f(x)$. (Gallian, Prout, & Winters, 1992) relaxed the definition of α-labelling and defined a weakly α-labelling to be a graceful labelling with the additional condition that there is a positive integer k such that for each edge xy, either $f(x) \leq k \leq f(y)$ or $f(y) \leq k \leq f(x)$. Researchers aimed to prove that all trees are graceful. This conjecture is known as the Ringel-Kotzig conjecture or the graceful tree conjecture (GTC) and despite the vast reservoir of papers that have been produced, it remains one of the most famous unsettled conjectures in graph theory.

The vast reservoir of papers that have been produced about graceful labelings, mostly exploit specific properties of some families of graphs as tools to prove that graphs belonging to these families are graceful. However, unfortunately, none of the techniques used so far are able to resolve the GTC. A list of families of trees that have been proved to be graceful can be found in almost every paper that has been written about the subject (See for example, (Luiz, Campos, & Richter, 2019; Rofa, 2021; Superdock, 2014; Wang, Yang, Hsu, & Cheng, 2015) and a detailed list of results is found in the dynamics survey Gallian (2018). As this paper is about rooted symmetric trees, we focus our literature review to this particular family of trees.

A rooted symmetric tree T is a rooted tree in which the degrees of all vertices at the same level or distance from the root vertex are equal. Figure 1 shows a rooted symmetric tree with root vertex v and seven levels including the root vertex in level 1. Rooted symmetric trees were repeatedly proved to be graceful by Bermond and Sotteau (Bermond & Sotteau, 1976) using G-designs, by Poljak and Sura (S. Poljak & Sˆura, 1982) and Robeva (Robeva, 2011) using algorithmic approaches and by Sandy et. al (Sandy, Rizal, Manurung, & Sugeng, 2018) using adjacency matrices and induction. However, all proofs are either inductive or algorithmic with many separate steps. In this paper we find a graceful labelling of rooted symmetric trees in the form of a direct algebraic function that *maps* each vertex of T to a unique label. The advantage of having such a direct functional proof lies in facilitating the production of graceful labellings of 'large' rooted symmetric trees in a computational way and in, relatively, minimal time. Furthermore, the graceful function that is found turns out to be a generalisation of the way a path is 'canonically' gracefully labelled as was algorithmically described by A. Rosa (Rosa, 1967). Notice that a path is a special type of a rooted symmetric tee. Therefore, from this point of view the result that is achieved in this paper can be considered as capturing Rosa's Algorithmic method for labelling paths in an algebraic way and generalising the result to include rooted symmetric trees. Figure 1(A) shows a graceful labelling of a 7-level multibranched rooted symmetric trees and Figure 1(B) shows the canonical algorithmic way of labelling a path (Which is a special case of a rooted symmetric tree) discovered by A. Rosa. In this paper we will show that both labellings can be derived through the similar algebraic functions.

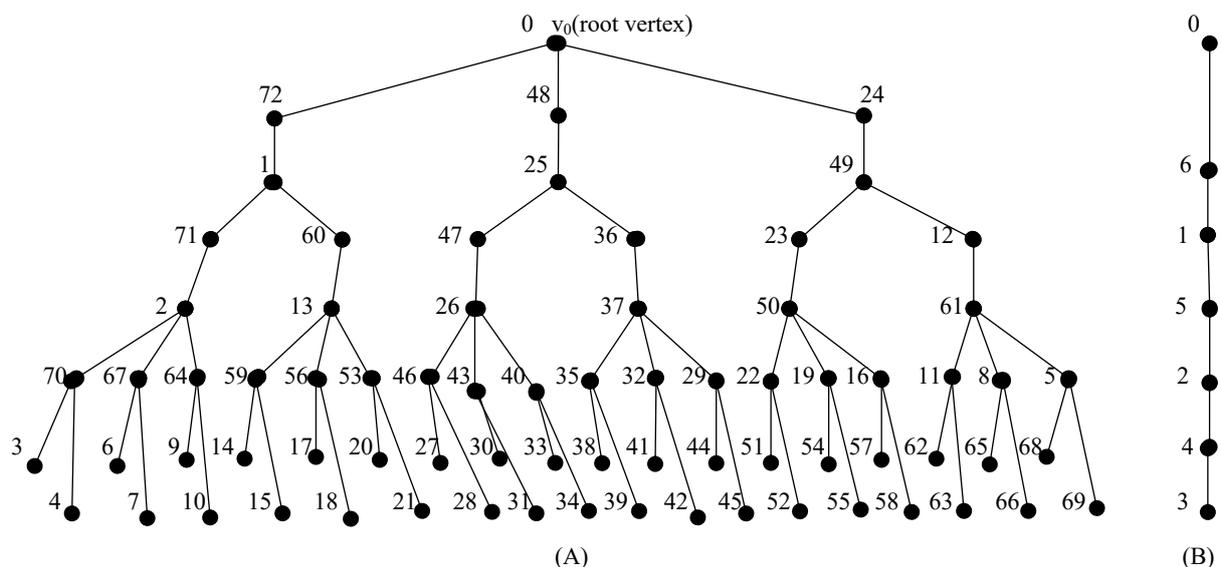

Figure 1

(A) Graceful labelling of a symmetric rooted tree with 7 levels (Including root vertex level).
(B) Canonical labelling of a path using the algorithmic method by A. Rosa.

## 2. Rooted Symmetric Trees: Definitions and Preliminary Results

In this section, we define the terminology and prove some preliminary results that will be required for providing an algebraic definition of a graceful labelling of a rooted symmetric tree (Section 3).

A rooted symmetric tree T is a rooted tree in which the degrees of all vertices at the same level are equal. Therefore, a rooted symmetric tree with q levels (level 1 to level q where level 1 contains root vertex v only) can be described by giving the constant degree at each level. For $1 \leq i \leq q$, let $k_i$ be the constant 'daughter' degree of each vertex at level i, that is, the number of vertices at level i+1 that are adjacent to a vertex at level i. Notice that vertices of the last level, that is level q, are not adjacent to any 'daughter' vertices, hence $k_q=0$. We call the sequence $S = (k_1, k_2, \ldots, k_{q-1})$ the daughter degree sequence of T and it, up to isomorphism totally identifies T. Label the $k_r$ daughter vertices of each vertex v at an arbitrary level, say level r, from left to right, by the numbers $0, 1, \ldots, k_r-1$ (not to be confused by the vertex labels of a graceful labelling of T). This labelling produces an embedding of T in which every vertex (except root vertex) is uniquely determined by a sequence of edge labels (not to be confused with the edge labels induced by a graceful labelling of T). For example, figure 2 shows a rooted symmetric tree with root vertex $v_0$ and degree sequence ($k_1 = 2, k_2 = 3, k_3 = 4$), in which vertices $v_2$, $v_5$, and $v_{32}$ are each uniquely represented by the edge sequences (1), (0,2), and (1,2,3), respectively. This identification of a vertex by a unique sequence of edge labels (the edges of the unique path that connects root vertex with the particular vertex that needs to be identified) is going to be used in algebraically defining the graceful labelling, f, of T.

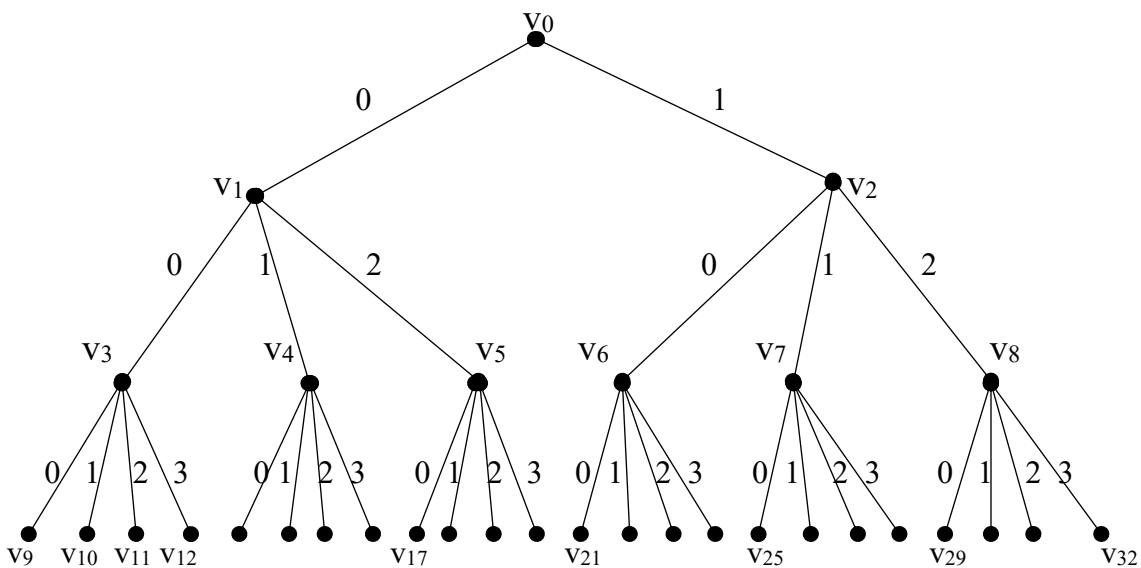

Figure 2. Identifying vertices by sequences

The next step is to find algebraic expressions for the number of vertices of the rooted, symmetric subtree having its root vertex at any given level, say level i between 1 and q, inclusive. We call this number *the $i^{th}$ level vertex size* and denote it by $h_i$. For example, in figure 1, $h_2$ denotes the 2nd level vertex size and is equal to the number of vertices of the rooted symmetric sub-tree with root vertex $v_2$ (or $v_1$). Using a simple counting technique, $h_2$ is counted as 1+3+3(4) = 16 (Remember, this tree has a daughter degree sequence ($k_1$=2, $k_2$=3, $k_3$=4). The following Lemma gives the level vertex size at any given level. The Lemma is easy to see by using the same counting technique that was used to find $h_2$ above.

**Lemma 1.** *For $1 \leq i \leq q-1$, we have* $h_i = 1 + k_i + k_i k_{i+1} + \ldots + k_i k_{i+1} \ldots k_{q-1}$ *and* $h_q = 1$

Notice that $h_1=|V|$ and $h_q$ is the number of vertices of the rooted symmetric tree having its root at the last level, which is a rooted symmetric tree consisting of the root vertex only, thus $h_q = 1$.

The following Corollary gives a recursive relation between two consecutive level vertex sizes

**Corollary 1**. For i between 1 and q-1, inclusive, we have, $h_{i+1} = \frac{h_i - 1}{k_i}$

Proof. Follows directly from Lemma 1. (It is useful for the reader to intuitively understand this formula by applying it to the rooted symmetric tree in figure 1).

Now, are in position to present an algebraic definition of a graceful labelling of a rooted symmetric tree. We do this in the following section.

## 3. Main Result: Algebraic Definition of a Graceful Labelling of Rooted Symmetric Trees.

As is the case throughout this paper, let T be a rooted symmetric tree with daughter degree sequence ($k_1$, $k_1$, ..., $k_{q-1}$) and root vertex $v_0$. Then using the set up that is described in the previous section, a vertex v at level r (except for the root vertex) is uniquely identified by the sequence of edge labels say ($x_1$, $x_2$, ..., $x_{r-1}$) of the edges that compose the unique path between the root vertex $v_0$ and vertex v. (See previous section for details and examples). We have the following main theorem.

**Theorem 1**. *T is graceful and the f: V $\rightarrow$ {0,1,..., |E|} defined by*
  1) $f(v_0) = 0$
  2) $f(v) = (k_1 - x_1)h_2 - x_2 h_3 - \ldots - x_{r-1} h_r - \frac{r-2}{2}$, *if r is even*
  3) $f(v) = x_1 h_2 + x_2 h_3 + \ldots + x_{r-1} h_r + \frac{r-1}{2}$, *if r is odd.*
*is a graceful labelling of the rooted symmetric tree T.*

Theorem 1 is proved by showing that the functions *f* and the function *g* : E → {1, ..., |E|} which is induced by f are both injective functions. but before doing that we demonstrate how Theorem 1 works in the following example. For the remainder of this paper f will be the function specified by theorem 1.

**Example 1**

Let T be the rooted symmetric tree in figure 1 with daughter degree sequence ($k_1$= 2, $k_2$ = 3, $k_3$ = 4). Then,
The level vertex sizes are:
$h_1 = 1 + k_1 + k_1 k_2 + k_1 k_2 k_3 = 33 = |V|$     [By Lemma 1].
$h_2 = 16$                                    [ By Lemma 1 or Corollary 2]
$h_3 = 5$                                     [ By Lemma 1 or Corollary 2]

$h_4 = 1$                                                             [ By Lemma 1 or Corollary 2]

Therefore, by the definition of f as per theorem 1, we have: (All these calculations can be performed computationally through using the algebraic function f of Theorem 1 in a simple code).

$f(v_0) = 0$

$v_1$ is represented by edge sequence $(x_1) = (0)$. Therefore, $f(v_1) = (k_1 - x_1)h_2 = 32$.
$v_2$ is represented by edge sequence $(x_1) = (1)$. Therefore, $f(v_2) = (k_1 - x_1)h_2 = 16$.
$v_3$ is represented by edge sequence $(x_1, x_2) = (0,0)$. Therefore, $f(v_3) = x_1 h_2 + x_2 h_3 + 1 = 1$.
$v_4$ is represented by edge sequence $(x_1, x_2) = (0,1)$. Therefore, $f(v_4) = x_1 h_2 + x_2 h_3 + 1 = 6$.
$v_5$ is represented by edge sequence $(x_1, x_2) = (0,2)$. Therefore, $f(v_5) = x_1 h_2 + x_2 h_3 + 1 = 11$.
$v_6$ is represented by edge sequence $(x_1, x_2) = (1,0)$. Therefore, $f(v_6) = x_1 h_2 + x_2 h_3 + 1 = 17$.
$v_7$ is represented by edge sequence $(x_1, x_2) = (1,1)$. Therefore, $f(v_7) = x_1 h_2 + x_2 h_3 + 1 = 22$.
$v_8$ is represented by edge sequence $(x_1, x_2) = (1,2)$. Therefore, $f(v_8) = x_1 h_2 + x_2 h_3 + 1 = 27$.
$v_9$ is represented by edge sequence $(0,0,0)$. Therefore, $f(v_9) = (k_1 - 0)h_2 - 0h_3 - 0h_4 - 1 = 31$.
$v_{10}$ is represented by edge sequence $(0,0,1)$. Therefore, $f(v_{10}) = (k_1 - 0)h_2 - 0h_3 - 1h_4 - 1 = 30$.
$v_{11}$ is represented by edge sequence $(0,0,2)$. Therefore, $f(v_{11}) = (k_1 - 0)h_2 - 0h_3 - 2h_4 - 1 = 29$.
$v_{12}$ is represented by edge sequence $(0,0,3)$. Therefore, $f(v_{12}) = (k_1 - 0)h_2 - 0h_3 - 3h_4 - 1 = 28$.
$v_{13}$ is represented by edge sequence $(0,1,0)$. Therefore, $f(v_{13}) = (k_1 - 0)h_2 - 1h_3 - 0h_4 - 1 = 26$.
$v_{14}$ is represented by edge sequence $(0,1,1)$. Therefore, $f(v_{14}) = (k_1 - 0)h_2 - 1h_3 - 1h_4 - 1 = 25$.
$v_{15}$ is represented by edge sequence $(0,1,2)$. Therefore, $f(v_{15}) = (k_1 - 0)h_2 - 1h_3 - 2h_4 - 1 = 24$.
$v_{16}$ is represented by edge sequence $(0,1,3)$. Therefore, $f(v_{16}) = (k_1 - 0)h_2 - 1h_3 - 3h_4 - 1 = 23$.
$v_{17}$ is represented by edge sequence $(0,2,0)$. Therefore, $f(v_{17}) = (k_1 - 0)h_2 - 2h_3 - 0h_4 - 1 = 21$.
$v_{18}$ is represented by edge sequence $(0,2,1)$. Therefore, $f(v_{18}) = (k_1 - 0)h_2 - 2h_3 - 1h_4 - 1 = 20$.
$v_{19}$ is represented by edge sequence $(0,2,2)$. Therefore, $f(v_{19}) = (k_1 - 0)h_2 - 2h_3 - 2h_4 - 1 = 19$.
$v_{20}$ is represented by edge sequence $(0,2,3)$. Therefore, $f(v_{20}) = (k_1 - 0)h_2 - 2h_3 - 3h_4 - 1 = 18$.
$v_{21}$ is represented by edge sequence $(1,0,0)$. Therefore, $f(v_{21}) = (k_1 - 1)h_2 - 0h_3 - 0h_4 - 1 = 15$.
$v_{22}$ is represented by edge sequence $(1,0,1)$. Therefore, $f(v_{22}) = (k_1 - 1)h_2 - 0h_3 - 1h_4 - 1 = 14$.
$v_{23}$ is represented by edge sequence $(1,0,2)$. Therefore, $f(v_{23}) = (k_1 - 1)h_2 - 0h_3 - 2h_4 - 1 = 13$.
$v_{24}$ is represented by edge sequence $(1,0,3)$. Therefore, $f(v_{24}) = (k_1 - 1)h_2 - 0h_3 - 3h_4 - 1 = 12$.
$v_{25}$ is represented by edge sequence $(1,1,0)$. Therefore, $f(v_{25}) = (k_1 - 1)h_2 - 1h_3 - 0h_4 - 1 = 10$.
$v_{26}$ is represented by edge sequence $(1,1,1)$. Therefore, $f(v_{26}) = (k_1 - 1)h_2 - 1h_3 - 1h_4 - 1 = 9$.
$v_{27}$ is represented by edge sequence $(1,1,2)$. Therefore, $f(v_{27}) = (k_1 - 1)h_2 - 1h_3 - 2h_4 - 1 = 8$.
$v_{28}$ is represented by edge sequence $(1,1,3)$. Therefore, $f(v_{28}) = (k_1 - 1)h_2 - 1h_3 - 3h_4 - 1 = 7$.
$v_{29}$ is represented by edge sequence $(1,2,0)$. Therefore, $f(v_{29}) = (k_1 - 1)h_2 - 2h_3 - 0h_4 - 1 = 5$.
$v_{30}$ is represented by edge sequence $(1,2,1)$. Therefore, $f(v_{30}) = (k_1 - 1)h_2 - 2h_3 - 1h_4 - 1 = 4$.
$v_{31}$ is represented by edge sequence $(1,2,2)$. Therefore, $f(v_{31}) = (k_1 - 1)h_2 - 2h_3 - 2h_4 - 1 = 3$.
$V_{32}$ is represented by edge sequence $(1,2,3)$. Therefore, $f(v_{32}) = (k_1 - 1)h_2 - 2h_3 - 3h_4 - 1 = 2$.

Figure 3 shows the graceful labelling that is obtained by applying theorem 1 to T. It can be readily verified that all vertex labels and induced edge labels are distinct.

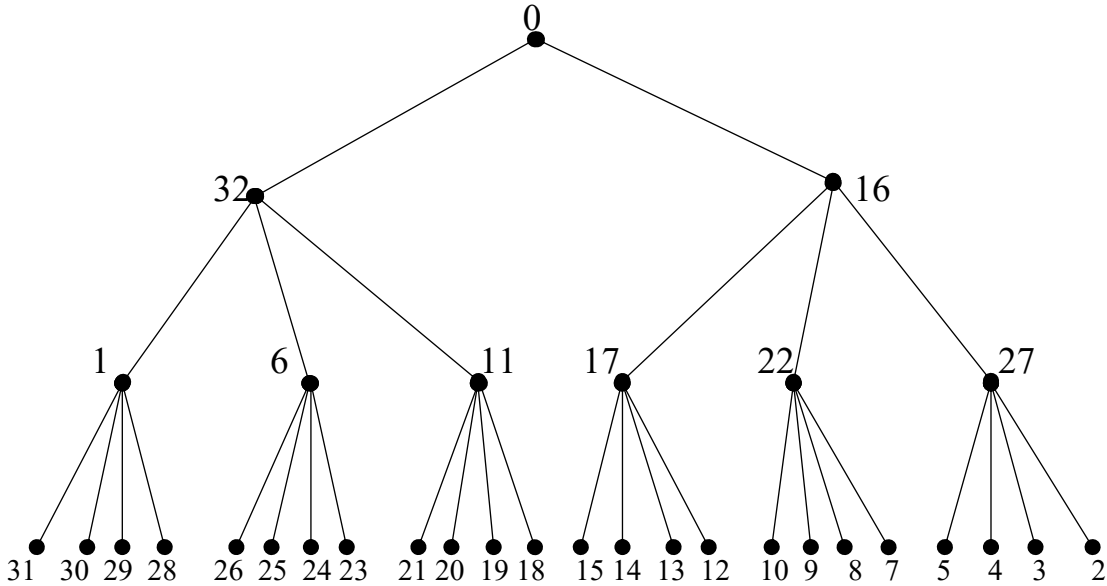

Figure 3: A graceful labelling of a rooted symmetric tree found by applying the *function* defined in Theorem 1.

**Proof of Theorem 1.**

Theorem 1 is proved by showing that functions *f: V → {0,1,…,|E|}*, as defined in Theorem 1, and *g: E → {1, … ,|E|}* where for edge uv g(uv)=|f(u)-f(v)|, are both injective functions. This is going to be achieved by presenting algorithms that show that the inverse of each of these functions are well defined functions.

Let f⁻¹: *{0,1,…,|E|} → VG* be the inverse *relation* of f. Our aim is to show that f is injective by showing that f⁻¹ is a well-defined function. That is, each element in the domain of f⁻¹ is mapped to exactly one vertex of T. So, let m∈*{0,1,…,|E|}*. The following algorithm will be shown to map m to unique vertex of T

**<u>Algorithm 1</u>**

Let m∈*{0,1,…,|E|}*.

(1) If m=0, then f⁻¹(m)= $v_0$, at level 1. End.
If m≠0, continue.

(2) <u>Testing whether f⁻¹(m)∈Level 2:</u>
Set $m' = k_1 h_2 - m$
Divide $m'$ by $h_2$. Let $\frac{m'}{h_2} = e_1$ $with\ remainder\ r_1'$.
If $r_1' = 0$, then, f⁻¹(m)∈Level 2, and is the vertex identified by the single-entry sequence ($e_1$). (Refer to section 2 and figure 1 to see how we associate vertices with sequences)
End.
*If* $r_1' \neq 0$. Continue.

(3) <u>Testing whether f⁻¹(m)∈Level 3</u>

Divide $m$ by $h_2$. Let $\frac{m}{h_2} = q_1$ $with$ $remainder$ $r_1 \neq 0$ (because if m is divisible $h_2$ so is $m'$, which is not the case by step 2)
Set $r_2 = r_1 - 1$.
Divide $r_2$ by $h_3$ and let $\frac{r_2}{h_3} = q_2$ with remainder $r_3$.
If $r_3 = 0$, then, f$^{-1}$(m)$\in$Level 3, and is the vertex identified by the sequence $(q_1, q_2)$. End.
$If$ $r_3 \neq 0$. Continue.

(4) <u>Testing whether f$^{-1}$(m)$\in$Level 4</u> (Algorithm continues from last step of (2))
Divide $r_1'$ by $h_3$ and let $\frac{r_1'}{h_3} = e_2$ with remainder $r_2' \neq 0$ (because $r_3 \neq 0$).
Set $r_3' = r_2' - 1$.
Divide $r_3'$ by $h_4$ and let $\frac{r_3'}{h_4} = e_3$ with remainder $r_4'$.
If $r_4' = 0$, then, f$^{-1}$(m)$\in$Level 4, and is the vertex identified by the sequence $(e_1, e_2, e_3)$. End.

$If$ $r_4' \neq 0$. Continue.

(5) <u>Testing whether f$^{-1}$(m)$\in$Level 5</u> (Algorithm continues from last step of (3))
Divide $r_3$ by $h_4$ and let $\frac{r_3}{h_4} = q_3$ with remainder $r_4 \neq 0$ (because $If$ $r_4' \neq 0$)
Set $r_5 = r_4 - 1$.
Divide $r_5$ by $h_5$ and let $\frac{r_5}{h_5} = q_4$ with remainder $r_6$.
If $r_6 = 0$, then, f$^{-1}$(m)$\in$Level 5, and is the vertex identified by the sequence $(q_1, q_2, q_3, q_4)$. End.
$If$ $r_6 \neq 0$. Continue.

(6) <u>Testing whether f$^{-1}$(m)$\in$Level 6</u> (Algorithm continues from last step of (4))
Divide $r_4'$ by $h_5$ and let $\frac{r_4'}{h_5} = e_4$ with remainder $r_5' \neq 0$ (because $r_6 \neq 0$).
Set $r_6' = r_5' - 1$
Divide $r_6'$ by $h_6$ and let $\frac{r_6'}{h_6} = e_5$ with remainder $r_7'$.
If $r_7' = 0$, then, f$^{-1}$(m)$\in$Level 6, and is the vertex identified by the sequence $(e_1, e_2, e_3, e_4, e_5)$. End.

$If$ $r_7' \neq 0$. Continue.

Repeat until a remainder of 0 is obtained and set f$^{-1}$(m) to be the unique vertex identified by the sequence ($e_1$, $e_2$, $e_3$, …) if the vertex belongs an even level and ($q_1$, $q_2$, $q_3$, …), if the vertex belongs to an odd level.
**End of Algorithm 1.**

Note that algorithm 1 assigns a unique vertex to each element in the domain of f$^{-1}$ and that it must terminate because in subsequent step, the remainder from the previous step is reduced by 1. This proves that f$^{-1}$ is a well-defined function and therefore, the function f is an injective function which is one of the first requirements for f to be a graceful labelling of T. A similar, but slightly more complicated algorithm, can show that g$^{-1}$ is a well-defined function. Therefore, f is a graceful labelling of T.

**End of Proof of Theorem 1.**

We illustrate Algorithm 1 in the following example.

**Example 2**
Consider the rooted symmetric tree T in figures 2 and the graceful labelling of T in figure 3.
a) Applying Algorithm 1 to vertex label m = 32 from figure 3 results in $f^{-1}(32) = (0)$ which correctly identifies the corresponding vertex $v_1$ in figure 2.
b) Applying Algorithm 1 to vertex label m = 10 from figure 3 results in $f^{-1}(10) = (1,1,0)$ which correctly identifies the corresponding vertex $v_{25}$ in figure 2.

Next, we will show that the canonical algorithmic graceful labelling method that was applied by A. Rosa to gracefully label paths can be captured by an algebraic function according to theorem 1. To clarify the concept, consider path P on 7 vertices, shown in figure 1(B). P can be thought of as a rooted symmetric tree with root vertex v (end vertex shown in the figure). By finding the level vertex sizes of P and applying the corresponding function f of theorem 1, we get the same labelling that could be found by using the algorithmic method of A. Rosa.

Our final result is to prove that the function f of theorem 1 is a weakly α-labelling of any rooted symmetric tree with daughter degree sequences $(k_1, k_2, …)$, where $k_1 = 2$

**Lemma 1.** Let T be a rooted symmetric tree with daughter degree sequence $(k_1, k_2, …)$, where $k_1 = 2$, then T has a weakly α-labelling.

**Proof.** Let T be a q level rooted symmetric tree with daughter degree sequence $(k_1, k_2, …)$, where $k_1 = 2$. By Theorem 1 the function f, as defined in the theorem is a graceful labelling of T. Let u and v be the adjacent vertices that are identified by the sequences (1) and (1,0), respectively. Then $f(u) = h_2$ and $f(v) = h_2 + 1$. Therefore, edge uv has edge label 1 and hence, the value of an integer that could satisfy a weakly alpha labelling is either $h_2$ or $h_2 + 1$. We claim that the integer $k = h_2$ satisfies the condition of f being a weakly alpha labelling.
**Proof of claim.** Without loss of generality, assume that vertex x is identified by the sequence $(x_1, x_2, ..., x_{r-1})$ where r≤q, and y is identified by the sequence $(x_1, x_2, ..., x_r)$. Also, without loss of generality, assume that r is odd. Then $f(x) = x_1 h_2 + x_2 h_3 + ... + x_{r-1} h_r + \frac{r-1}{2}$ and $f(y) = (2 - x_1) h_2 - x_2 h_3 - ... - x_{r-1} h_r - x_r h_{r+1} - \frac{r-1}{2}$. Now $x_1$ is either 0 or 1 (because $k_1 = 2$). If $x_1 = 0$ then $f(x) \leq h_2$ because, by using a counting principle we have $x_2 h_3 + ... + x_{r-1} h_r + \frac{r-1}{2} \leq h_2$
(this could also be proved by induction), and for a similar reason, $h_2 \leq f(y)$. If $x_1 = 1$, then it can similarly be shown that $f(y) \leq h_2 \leq f(x)$. Therefore, in both cases either $f(x) \leq h_2 \leq f(y)$ or $f(y) \leq h_2 \leq f(x)$ and hence, f is a weakly α-labelling of T.
End of proof.

**Corollary 1.** Binary trees have weakly α-labelling.
**Proof.** The proof follows from corollary 1 as a binary tree is a rooted symmetric tree with daughter degree sequence (2, 2, …).

# References

Bermond, J. C., & Sotteau, D. (1976). Graph decompositions and G-designs. *In 5th British Combinatorial Conference, 1975, Congressus Numerantium 15, Utilitas Math Pub*, 53-72.

Gallian, J. A. (2018). A dynamic survey of graph labeling. *Electron. J. Combin. DS6* 1-502.

Gallian, J. A., Prout, J., & Winters, S. (1992). Graceful and harmonious labelings of prisms and related graphs. *Ars Combin., 34*, 213-222.

Golomb, S. W. (1972). How to Number a Graph. *Graph Theory and Computing, R. C. Read, ed., Academic Press, New York*, 23-37.

Luiz, A. G., Campos, C. N., & Richter, R. B. (2019). On $\alpha$-labellings of lobsters and trees with a perfect matching. *Discrete Applied Mathematics, 268*, 137-151. doi:10.1016/j.dam.2019.05.004

Ringel, G. (1963). Problem 25. *In: Theory of Graphs and its Applications. Proc. Symposium Smolenice Prague (1964), 162*.

Robeva, E. (2011). An extensive survey of graceful trees. *Udergraduate Honours Thesis, Standford University, USA*.

Rofa, R. (2021). Graceful and Strongly Graceful Permutations. . Retrieved from https://arxiv.org/pdf/2107.14007.pdf

Rosa, A. (1967). On certain valuations of the vertices of a graph. *In: Theory of Graphs. Proceedings of the Symposium, 1966, (Rome), Gordon and Breach, New York*, 349–355.

S. Poljak, S., & Sˆura, M. (1982). An algorithm for graceful labeling of a class of symmetrical trees. *Ars Combin., 14*, 57-66.

Sandy, I. P., Rizal, A., Manurung, E. N., & Sugeng, K. A. (2018). Alternative construction of graceful symmetric trees. *Journal of Physics.: Conf. Ser. 1008 012031*. doi:10.1088/1742-6596/1008/1/012031

Superdock, M. C. (2014). The graceful tree conjecture: A class of graceful diameter-6 trees. *Ithaca: Cornell University Library, arXiv.org*.

Wang, T.-M., Yang, C.-C., Hsu, L.-H., & Cheng, E. (2015). Infinitely many equivalent versions of the graceful tree conjecture. *Applicable Analysis and Discrete Mathematics, 9*(1), 1-12. doi:10.2298/aadm141009017w